\documentclass[11pt]{article}

\usepackage{enumerate}
\usepackage{amsmath,xspace,amssymb,mathrsfs,theorem}

\input xy
\xyoption{all}
\xyoption{2cell}
\UseAllTwocells
\CompileMatrices

\title{Notions of Lawvere theory \thanks{Both authors gratefully acknowledge the
support of the Australian Research Council; the second-named author also gratefully acknowledges the support of the Ministry of Education of the Czech Republic by the project MSM 00216224409 and of the Grant Agency of the Czech
Republic by the project 201/06/0664.} }
\author{
Stephen Lack \\
School of Computing and Mathematics \\
University of Western Sydney \\
Locked Bag 1797 Penrith South DC NSW 1797\\
Australia\\
{\tt s.lack@uws.edu.au}
\and 
Ji\v r\'i Rosick\'y \\
Department of Mathematics and Statistics \\
Masaryk University
Kotl\'a\v rsk\'a 2 60000 Brno\\
Czech Republic\\
{\tt rosicky@math.muni.cz}
}

\renewcommand{\phi}{\varphi}

\newcommand{\A}{{\ensuremath{\mathscr{A}}}\xspace}
\newcommand{\B}{{\ensuremath{\mathscr{B}}}\xspace}
\newcommand{\C}{{\ensuremath{\mathscr{C}}}\xspace}
\newcommand{\D}{{\ensuremath{\mathscr{D}}}\xspace}
\newcommand{\DD}{{\ensuremath{\mathbb{D}}}\xspace}
\newcommand{\E}{{\ensuremath{\mathscr{E}}}\xspace}
\newcommand{\F}{{\ensuremath{\mathscr{F}}}\xspace}
\newcommand{\G}{{\ensuremath{\mathscr{G}}}\xspace}
\newcommand{\GG}{{\ensuremath{\overline{\mathscr{G}}}}\xspace}
\newcommand{\I}{{\ensuremath{\mathscr{I}}}\xspace}
\newcommand{\K}{{\ensuremath{\mathscr{K}}}\xspace}

\newcommand{\Ksf}{{\ensuremath{\mathscr{K}_{sf}}}\xspace}
\newcommand{\LL}{{\ensuremath{\mathscr{L}}}\xspace}
\newcommand{\M}{{\ensuremath{\mathscr{M}}}\xspace}
\renewcommand{\P}{{\ensuremath{\mathscr{P}}}\xspace}

\renewcommand{\S}{{\ensuremath{\mathscr{S}}}\xspace}
\newcommand{\T}{{\ensuremath{\mathscr{T}}}\xspace}
\newcommand{\FT}{{\ensuremath{\mathscr{F\!T}}}\xspace}

\newcommand{\V}{{\ensuremath{\mathscr{V}}}\xspace}
\newcommand{\Vf}{{\ensuremath{\mathscr{V}_f}}\xspace}
\newcommand{\Vsf}{{\ensuremath{\mathscr{V}_{sf}}}\xspace}

\newcommand{\II}{{\ensuremath{\mathbb{I}}}\xspace}

\newcommand{\PA}{\ensuremath{\mathscr{P\!A}}\xspace}
\newcommand{\PK}{\ensuremath{\mathscr{P\!K}}\xspace}

\newcommand{\FPP}{\textnormal{\bf FPP}\xspace}
\newcommand{\Mod}{\textnormal{\bf Mod}\xspace}
\newcommand{\PhiMod}{\textnormal{\bf $\Phi$-Mod}\xspace}
\newcommand{\PhiCts}{\textnormal{\bf $\Phi$-Cts}\xspace}
\newcommand{\Lex}{\textnormal{\bf Lex}\xspace}
\newcommand{\Law}{\textnormal{\bf Law}\xspace}
\newcommand{\Mnd}{\textnormal{\bf Mnd}\xspace}
\newcommand{\mnd}{\textnormal{\bf mnd}\xspace}
\renewcommand{\th}{\textnormal{\bf th}\xspace}

\newcommand{\Set}{\textnormal{\bf Set}\xspace}
\newcommand{\Setf}{\ensuremath{\textnormal{\bf Set}_{\textnormal{\bf f}}}\xspace}
\newcommand{\SSet}{\textnormal{\bf SSet}\xspace}
\newcommand{\Ab}{\textnormal{\bf Ab}\xspace}

\newcommand{\Gpd}{\textnormal{\bf Gpd}\xspace}
\newcommand{\Gph}{\textnormal{\bf Gph}\xspace}
\newcommand{\RGph}{\textnormal{\bf RGph}\xspace}
\newcommand{\Cat}{\textnormal{\bf Cat}\xspace}
\newcommand{\CGTop}{\textnormal{\bf CGTop}\xspace}
\newcommand{\VCat}{\textnormal{\bf \V-Cat}\xspace}

\newcommand{\op}{\ensuremath{^{\textnormal{op}}}}
\newcommand{\ev}{\textnormal{ev}}
\newcommand{\Fam}{\textnormal{Fam}}
\newcommand{\Lan}{\textnormal{Lan}}
\newcommand{\colim}{\textnormal{colim}}

\newcommand{\ot}{\otimes}

\newdir{ >>}{{}*!/-10pt/@{>>}}
\newdir{ >}{{}*!/-8pt/@{>}}

\renewcommand{\t}{\times}         
\def\endproof{{\parfillskip=0pt\hfill$\Box$\vskip 10pt}}

\newcommand{\two}{\ensuremath{{\hbox{\textrm 2}\kern-.25em
        \hbox{\vrule height1.5ex width 0.4pt depth -.2ex}}\kern.2em}\xspace}

\newtheorem{theorem}{Theorem}[section]    
\newtheorem{corollary}[theorem]{Corollary}   
\newtheorem{proposition}[theorem]{Proposition}   
\newtheorem{lemma}[theorem]{Lemma}   
\newtheorem{definition}[theorem]{Definition}   

\newtheorem{axiom}{Axiom}
\renewcommand{\theaxiom}{\Alph{axiom}.}

{\theorembodyfont{\rmfamily}
\newtheorem{remark}[theorem]{Remark}   
\newtheorem{example}[theorem]{Example}   
}

\newcommand{\proof}{\noindent{\sc Proof:}\xspace}

\newcommand{\lsfp}{{locally strongly finitely presentable}\xspace}
\newcommand{\sfp}{{strongly finitely presentable}\xspace}
\newcommand{\sfly}{{strongly finitely}\xspace}
\newcommand{\sfinite}{{strongly finite}\xspace}
\newcommand{\Sfinite}{{Strongly finite}\xspace}

\begin{document}

\label{firstpage}
\maketitle

\begin{abstract}
Categorical universal algebra can be developed either using Lawvere 
theories (single-sorted finite product theories) or using monads, and the category of Lawvere theories is
equivalent to the category of finitary monads on \Set. We show how 
this equivalence, and the basic results of universal algebra, can be
generalized in three ways: replacing \Set by another category,
working in an enriched setting, and by working with another class of 
limits than finite products. 
\end{abstract}


Classical universal algebra begins with structures of the following 
type: a set $X$ equipped with operations $X^n\to X$ for various 
natural numbers $n$, subject to equations between induced operations. 
For example the structure of group can be encoded using a binary 
operation $m:X^2\to X$ (multiplication), a unary operation $(~)^{-1}:X\to X$ 
(inverse) and a nullary one $i:1\to X$ (unit). Then there is an equation
encoding associativity, two equations for the two unit laws, and two equations
for the inverses.

There are two main ways for treating such structures categorically. In 
each case, one ends up giving all the operations generated, in a suitable 
sense,  by a presentation as in the previous paragraph. Thus the structure
of group (as opposed to some particular group) becomes a mathematical object
in its own right.

On the one hand, for each such type of structure there is a small category \T 
with finite products, such that a particular instance of the structure 
is a finite-product-preserving functor from \T to \Set. The objects of 
\T will be exactly the powers $T^n$ of a fixed object $T$. Such a category 
is called a {\em Lawvere theory}. In this case, the operations in the 
structure are seen as the morphism in the theory \T.

On the other hand, for each such type of structure there is a monad $T$
on \Set, such that a particular instance of the structure is an algebra
for the monad. The basic idea this time is to describe the structure in 
question by saying what the free algebras are. Since we have been considering
only {\em finitary} operations (operations $X^n\to X$ where $n$ is finite)
the resulting monad $T$ will also be finitary. This means that it is 
determined by its behaviour on finite sets; formally, this means
that the endofunctor $T$ preserves filtered colimits.

There is a very well developed theory of Lawvere theories and of finitary
monads, and of the equivalence between the two notions. In this paper
we are interested in what happens when we move beyond the context of 
structures borne by a set, to consider structures borne by other objects. 
Lawvere theories are single-sorted theories; much of the time we shall
spend on the 
many-sorted case, or rather the ``unsorted'' case, where a theory consists
of an arbitrary small category with finite products, and a model is just
a finite-product-preserving functor.

We shall consider three different ways in which the classical setting
could be adapted; the three ways are to some extent independent of
each other.

\paragraph{Replace the base category \Set by some other suitable category \K}

We could consider structures borne not by sets but by objects of other categories. For 
example we could take \K to be the category \Gph of (directed) graphs, and 
consider structures borne by graphs. An obvious example is a category:
this is a graph equipped with a composition law and identities.

From the point of view of monads, this change is straightforward enough: one
can simply consider monads on \Gph rather than monads on \Set. As for theories,
there is no problem considering models of a Lawvere theory \T in \Gph, one
simply takes finite-product-preserving functors from \T to \Gph rather than 
from \T to \Set. But often this is too restrictive; for example it would not 
include the case of categories. One needs to allow more general arities than
natural numbers; in fact it is (certain) objects of \Gph which will serve
as arities.

\paragraph{Deal not with ordinary categories but with categories enriched over some
suitable monoidal category \V}

A monoidal category can clearly be thought of as a category with extra 
structure, of an algebraic sort. When we come to think of the morphisms
between such structures, this naive algebraic point of view would lead us
to take as morphisms the strict monoidal functors (those which preserve the
unit and tensor product in the strict sense of equality: $F(X\ot Y)=FX\ot FY$).
Although such strict monoidal functors have an important role to play, much
more important are various ``weak'' notions, where the monoidal structure
is preserved up to coherent isomorphism, or perhaps only up to a coherent
non-invertible comparison map. These weaker notions of homomorphism can be 
recovered by working not over the category of small categories and functors,
but over the 2-category of small categories, functors, and natural 
transformations. In other words, we work with categories enriched in \Cat.

More generally one could consider many other monoidal categories \V, such
as abelian groups, chain complexes, simplicial sets, and many more:
see \cite{Kelly-book}. 

When we are working over \V, the ``default'' choice of \K becomes \V 
rather than \Set, but many other \V-categories may also be suitable.

\paragraph{Replace finite products by some suitable class $\Phi$ of limits.}

This is probably the most important aspect for this paper.

Some structures cannot be defined with finite products alone. Once again
the structure of categories is a good example. An internal category consists
of an object-of-objects $C_0$, an object-of-morphisms $C_1$, with domain 
and codomain maps $d,c:C_1\to C_0$, together with identities $i:C_0\to C_1$
and composition $m:C_2\to C_1$ satisfying various equations, where $C_2$
is object-of-composable-pairs, given by the pullback
$$\xymatrix{
C_2 \ar[r] \ar[d] & C_1 \ar[d]^{c} \\ C_1 \ar[r]_{d} & C_0. }$$
So to define internal categories one needs a class $\Phi$ of limits
which contains pullbacks.

Once again, the limits in $\Phi$ should be weighted limits in the 
\V-enriched sense. 

Our main new example of such a class is where $\Phi$ consists of 
what we call \sfinite limits. This is a class intermediate between
finite products and finite limits. An object $X\in\V$ is
{\em \sfinite} if the internal hom functor
$[X,-]:\V\to\V$ preserves sifted colimits (see \cite{sifted} or
Section~\ref{sect:fp} below). The \sfinite limits are those
which can be constructed using finite products and powers (cotensors)
by \sfinite objects. These correspond to monads which
preserve sifted colimits. 
As observed for example in \cite{Rezk-thesis}, every monad on a cocomplete
symmetric monoidal closed category arising from an operad preserves 
filtered colimits and reflexive coequalizers; thus, using 
Proposition~\ref{prop:funny}, it preserves sifted colimits.

It is worth observing that as well as being the free completion of 
the category \Setf of finite sets under filtered colimits, \Set is also 
the free completion of \Setf under sifted colimits. Thus an endofunctor
of \Set preserves filtered colimits if and only if it preserves sifted
colimits. Thus the well-known equivalence between Lawvere theories and
finitary monads on \Set could also be seen as an equivalence between
Lawvere theories and sifted-colimit-preserving monads on \Set. This is 
a very special property of \Set; for other base categories, one therefore
potentially has (at least) two ways to generalize the equivalence. 

\paragraph{}

So altogether we have three changes
\begin{enumerate}[(A)]
\item replace the base category \Set by some other suitable category \K
\item deal not with ordinary categories but with categories enriched over
  some suitable monoidal category \V
\item replace finite products by some suitable class $\Phi$ of limits
\end{enumerate}
and these appear to be independent to some extent, apart from the obvious
restrictions that \K should be a \V-category and the limits in $\Phi$
should be weighted \V-limits. But the real relationship between these
conditions appears when we try to say what ``suitable'' should mean in
each case; doing this will be one of the main aims of the paper. 

We'll
see that a key aspect is that it gives a sort of decomposition or factorization
of all limits into two parts. In fact it's easier to think in terms of colimits. For example if $\Phi$ consists of the finite limits, then 
we have the decomposition of arbitrary colimits as filtered colimits of
finite ones. Then again, if $\Phi$ consists of the finite products we have
a corresponding decomposition into sifted colimits of finite coproducts
 (see Section~\ref{sect:sifted} for more about sifted colimits).

There is a trade-off between different levels of generality of theory. 
As usual, the more general the notion of theory, the more general the categories
of models that can be described, but the less that can be proved about them.
There are also more subtle effects. For example, consider the case $\V=\Set$ 
of unenriched categories and compare finite limit theories with finite product
theories. Finite limit theories are of course more general. On the other 
hand, in the first instance one needs a category with finite limits as
the category \K in which models are taken; whereas for a finite product 
theory, \K need only have finite products. As one goes further, however, 
it turns out that for the strongest results, \K should itself be the category
of models for a theory of the given type, and this is a stronger condition
on \K when we use finite products than it is when we use finite limits.

The main results of universal algebra that we
should like to obtain in our general setting are:

\paragraph{$\Phi$-algebraic functors have left adjoints.}

If $G:\S\to\T$ is morphism of theories, then there is an induced 
map $G^*:\Mod(\T)\to\Mod(\S)$, and such a $G^*$ we call $\Phi$-algebraic
(or just algebraic if $\Phi$ is understood). Our first basic fact
is that such a $G^*$ has a left adjoint; or rather, we give an explicit
construction of a left adjoint in terms of colimits in the base category
\K --- the existence of the left adjoint is known for general reasons.

This includes the existence (and construction) of free algebras for 
single-sorted theories.

\paragraph{The reflectiveness of models.}

For any theory \T, the category $\Mod(\T,\K)$ of algebras in \K is reflective 
in $[\T,\K]$.
This means that $\Mod(\T,\K)$ will be complete and cocomplete provided
that \K is so, as we shall usually suppose to be the case.
But we shall do more than just prove the reflectivity, we shall 
construct a reflection in terms of colimits in \K, and so obtain
a description of colimits in $\Mod(\T,\K)$ in terms of colimits in \K.

In the final chapter of \cite{Kelly-book}, the reflectiveness was 
proved under extremely general conditions, but without an explicit 
construction. Our setting will be much more restrictive, but will
allow an explicit construction (as explicit as is our knowledge of 
colimits in \K).

Although we can obtain an explicit construction of all colimits in 
$\Mod(\T,\K)$, it is particularly simple for the class of colimits which
commute in \K with $\Phi$-limits, since these colimits are formed in 
$\Mod(\T,\K)$ as in $[\T,\K]$, and so no reflection is required.

Once again we see the trade-off: the smaller the class of limits in $\Phi$,
the closer $\Mod(\T,\K)$ is to $[\T,\K]$, and so the greater our knowledge
of colimits in $\Mod(\T,\K)$. 

For example, if we take $\K=\V=\Set$ and compare the case where $\Phi$ consists
of all finite limits with that where it consists only of finite products,
in the first case the inclusion $\Mod(\T,\Set)\to[\T,\Set]$ preserves 
filtered colimits, while in the second case it also preserves reflexive
coequalizers.

It is perhaps worth noting that the reflectivity of models is a special
case of the existence of left adjoints to algebraic functors. Let $\FT$
be obtained by freely adding $\Phi$-limits to \T. Since \T already has
$\Phi$-limits, the canonical inclusion $\T\to\FT$ has a right adjoint $R$,
and now the induced algebraic functor $R^*:\Mod(\T,\K)\to\Mod(\FT,\K)$
may be identified with the inclusion $\Mod(\T,\K)\to[\T,\K]$, and so the 
left adjoint to $R^*$ gives the desired reflection.

On the other hand, we also have a sort of converse: the left adjoint 
to $G^*:\Mod(\T)\to\Mod(\S)$ can be obtained by first taking the left
Kan extension along $G$ and then reflecting into the subcategory of models.

\paragraph{The correspondence between (single-sorted) theories and monads.}

For suitable \V and $\Phi$, we consider the class of all colimits which
commute in \V with $\Phi$-limits; these will be called {\em $\Phi$-flat}.
We shall say that a \V-functor is {\em $\Phi$-accessible} if it preserves
$\Phi$-flat colimits, and that a \V-monad is $\Phi$-accessible if its 
underlying \V-functor is so.

For suitable \V-categories \K, we shall describe a notion of {\em \K-based 
$\Phi$-theory}, or {\em $\Phi$-theory in \K},  and prove that the category of these is equivalent to the 
category of $\Phi$-accessible \V-monads on \K, and furthermore that the 
algebras for a $\Phi$-accessible monad are the same as the models for the
corresponding $\Phi$-theory. This generalizes the classical 
equivalence between Lawvere theories and finitary monads on \Set, and is 
one of the main results of the paper.




\subsection*{Outline of paper}

We begin in Section~\ref{sect:enriched} with a review of the  enriched
category theory that will be needed in the paper. In Section~\ref{sect:axioms}
we describe the basic assumptions we shall make connecting our monoidal
category \V, our base \V-category \K, and the class of limits $\Phi$ to 
be considered. Section~\ref{sect:sifted} is largely a review of various ideas
relating to sifted colimits: these are the colimits which commute in \Set with
finite products. In Section~\ref{sect:examples} we describe the various possible
classes of limits, and the corresponding requirements on \V (and on \K). 
The last two sections contain our main results; 
we have divided these into those which are independent of the sorts, in 
Section~\ref{sect:general}, and those which relate specifically to single-sorted
theories, in Section~\ref{sect:Lawvere}.

\section{Review of enriched category theory}\label{sect:enriched}

We work over a symmetric monoidal closed category $\V=(\V_0,\ot,I)$
whose underlying ordinary category $\V_0$ is complete and cocomplete.
The general results in \cite{Kelly-book} on reflectivity of models, referred
to above, used the further assumption that \V is locally bounded, in the
sense of \cite[Chapter~6]{Kelly-book}. This includes all the key examples
of \cite{Kelly-book}, including the categories of sets, pointed sets,
abelian groups, 
modules over a commutative ring, chain complexes, categories, groupoids, 
simplicial sets, compactly generated  spaces (Hausdorff or otherwise, pointed 
or otherwise), Banach spaces, sheaves on a site, truth values,
and Lawvere's poset of extended non-negative real numbers. We shall often
make stronger assumptions on \V.

If the base \V is clear, we generally omit the prefix ``\V-'' and speak 
simply of a category, functor, or natural transformation. 

A {\em weight} will be a presheaf $F:\A\op\to\V$ which is a small
colimit of representables. If \A is small, any presheaf on \A is a small
colimit of representables, so there is no restriction. We sometimes say
that $F$ is small to mean that it is a small colimit of representables.
See \cite{small} for more on small functors.

For 
a \V-functor $S:\A\op\to\K$, the limit $\{F,S\}$ of $S$ weighted by $F$
is defined by a natural isomorphism
$$\K(X,\{F,S\})\cong[\A\op,\V](F,\K(X,S))$$
while for a \V-functor $R:\A\to\K$, the colimit $F*R$ of $R$
weighted by $F$ is defined by a natural isomorphism
$$\K(F*R,X)\cong[\A\op,\V](F,\K(R,X)).$$
Note that the previous two displayed equations appear to involve hom-objects
in $[\A\op,\V]$. If \A is large then $[\A\op,\V]$ does not exist as a 
\V-category; nonetheless, the desired hom-object will exist as an object
of \V, since $F$ is small (see \cite{small} for example).

In particular, if $F$ and $G$ are both small, then $[\A\op,\V](F,G)$
exists as an object of \V, and so we do have a \V-category $\PA$ of all 
small presheaves on \A. This is the free completion of \A under colimits
\cite{lindner-morita}.

A class $\Phi$ of limits means a class of weights; then 
$\Phi$-completeness or $\Phi$-continuity means the existence or the existence
and preservation of all limits with weights in $\Phi$. 

An important special case is the weight $C:\I\op\to\V$, where \I is
the unit \V-category consisting of a single object $*$ with hom
$\I(*,*)=I$. Then to give the weight is just to give an object $C\in\V$.
The $C$-weighted limit of a \V-functor $S:\I\op\to\K$ (that is, of an object
$S\in\K$ is defined by a natural isomorphism
$$\K(X,C\pitchfork S)\cong \V(C,[X,S])$$
has traditionally been called a {\em cotensor}, but we shall simply call
a {\em power}, or $C$-power where necessary. The corresponding colimit,
written $C\cdot S$, used to be called a {\em tensor}, but we shall call a
{\em copower}.

The ordinary, unweighted notion of limit can be seen as a special case.
Let \D be an ordinary category, and let $\F\D$ be the free \V-category on
\D. Then \V-functors $\F\D\to\K$ correspond to functors $\D\to\K_0$, and
we define the limit {\em in \K} of a functor $S:\D\to\K_0$ to be the 
limit of 
the corresponding $R:\F\D\to\K$ weighted by the terminal weight 
$\Delta 1:\F\D\to\K$. Limits of this form are called {\em conical}.

The universal property of $\{\Delta 1,R\}$ involves
an isomorphism in \V, and 
is strictly stronger than the universal property  of $\lim S$, which involves
only a bijection of sets. Nonetheless, the universal property of $\lim S$ 
does serve to identify $\{\Delta 1,R\}$ if the latter is known to exist. Thus
if \K and \LL have the relevant
limits, then to say that $F:\K\to\LL$ {\em preserves} a particular conical
limit is equivalent to saying that $F_0:\K_0\to\LL_0$ does so.

We shall need the following basic result:

\begin{proposition}\label{prop:Lan-weight}
Let $F:\A\op\to\V$ be a weight, and $J:\A\to\B$ and $S:\B\to\C$ functors.
Then 
$$F*SJ\cong \Lan_JF*S$$
either side existing if the other does.   
\end{proposition}

\proof
Here $\Lan_J F$ denotes the left Kan extension of $F:\A\op\to\V$ along
$J:\A\to\B\op$. 
The result follows from the calculation
\begin{align*}
\C(F*SJ,C) &\cong [\A\op,\V](F,\C(SJ,C)) \\
              &\cong [\B\op,\V](\Lan_JF,\C(S,C)) \\
              &\cong \C(\Lan_JF*S,C).
\end{align*}
\endproof

\begin{definition}
A class $\Phi$ of weights is said to be {\em saturated} if for any 
diagram $S:\D\to[\A\op,\V]$
in a presheaf category, with each $SD:\A\op\to\V$ in $\Phi$, and for 
any $F:\D\op\to\V$ in $\Phi$, the colimit $F*S\in[\A\op,\V]$ is also
in $\Phi$.
\end{definition}

The original reference \cite{Albert-Kelly} used the word ``closed'' in 
place of ``saturated'', but the latter is now standard. The basic result 
about a saturated class $\Phi$ is that the full subcategory of $[\A\op,\V]$ 
consisting
of the presheaves in $\Phi$ is the free completion of \A under 
$\Phi$-colimits (provided that \A is small).

A \V-functor $F:\A\to\B$ with small domain is said to be {\em dense}
\cite[Chapter~5]{Kelly-book}
if the induced \V-functor $\B(F,1):\B\to[\A\op,\V]$ sending 
$B\in\B$ to $\B(F-,B):\A\op\to\V$ is fully faithful.

Let $F:\A\op\to\V$ and $G:\B\op\to\V$ be weights, and \K a category
with $F$-limits and $G$-colimits. For any $S:\A\op\ot\B\to\K$,
and any $B\in\B$, we can form the limit $\{F,S(-,B)\}$ in \K, and this
defines the object part of a functor 
$\{F,S\}:\B\to\K:B\mapsto\{F,S(-,B)\}$, to which we
can now apply $G*-$ to obtain $G*\{F,S\}\in\K$. Similar we can
form a functor $G*S:\A\op\to\K$ and then $\{F,G*S\}\in\K$, and there
is a canonical comparison map
$$G*\{F,S\}\to\{F,G*S\}$$
in \K. If this is invertible for all $S$, we say that $F$-limits commute
in \K with $G$-colimits. The following observation was 
made in \cite{Kelly-Schmitt} in the case $\K=\V$:

\begin{proposition}
Let the \V-category \K have all $F$-limits  and all
$G$-colimits. Then the following are equivalent:
\begin{enumerate}
\item $F$-limits commute in \K with $G$-colimits;
\item $\{F,-\}:[\A,\K]\to\K$ is $G$-cocontinuous;
\item $G*-:[\B\op,\K]\to\K$ is $F$-continuous.
\end{enumerate}
\end{proposition}

More generally, if $\Phi$ and $\Psi$ are classes
of weights, we say that $\Phi$-limits commute in \K with $\Psi$-colimits
if this is so for all $F\in\Phi$ and all $G\in\Psi$. 

If $G$ commutes in \V with $\Phi$-limits we follow \cite{Kelly-Schmitt}
in calling  $G$ is 
{\em $\Phi$-flat}, or just $F$-flat if $\Phi=\{F\}$. This is by 
analogy with the case where $\V=\Ab$ and $\Phi$ consists of the finite 
conical limits. Then a one-object \Ab-category $B$ is a ring, and a 
weight $G:B\op\to\Ab$ is a $B$-module, while $G*-$ corresponds to 
tensoring over $B$. A module is flat exactly when tensoring with that 
module preserves finite limits.

In \cite[Proposition~5.4]{Kelly-Schmitt}, the class of $\Phi$-flat weights
is shown to be saturated. An important part of this is the following:

\begin{proposition}\label{prop:flat-saturated}
If $F:\A\op\to\V$ if $\Phi$-flat and $G:\A\op\to\B\op$ is 
arbitrary, then $\Lan_GF:\B\op\to\V$ is $\Phi$-flat.
\end{proposition}

\proof
To say that $\Lan_G F$ is $\Phi$-flat is to say that $\Lan_G F*-$ 
is $\Phi$-continuous. But by Proposition~\ref{prop:Lan-weight}, 
this $\Lan_G F*-$ is given by the composite
$$\xymatrix{
[\B,\V] \ar[r]^{[J,\V]} & [\A,\V] \ar[r]^{F*-} & \V }$$
and $[J,\V]$ preserves all limits, since limits in presheaf categories are
calculated pointwise, while $F*-$ is $\Phi$-continuous since $F$
is $\Phi$-flat.
\endproof

\section{Key requirements}\label{sect:axioms}


The first requirement involves $\Phi$ and \V. It is convenient to suppose
that $\Phi$ is {\em locally small} \cite{Kelly-Schmitt} in the sense that
for any small \A, the closure of the representables in $[\A\op,\V]$ under
$\Phi$-colimits is again small; typically this will happen because $\V_0$ is
locally presentable and all
weights in $\Phi$ are $\alpha$-presentable for some regular cardinal 
$\alpha$.

\begin{axiom} If \A is a small \V-category with $\Phi$-limits, and
$F:\A\to\V$ is $\Phi$-continuous, then so is $F*-:[\A\op,\V]\to\V$.
\end{axiom}

Note that $F*-$ is the left Kan extension $\Lan_YF$ of $F$ along
the Yoneda embedding.
This condition has been considered by many different authors in various
special cases, and some of these are listed below when we turn to examples.
In particular, it holds in the case $\V=\Set$ if $\Phi$ consists of either
finite products or finite limits. It was considered, still in the 
case $\V=\Set$, for a general class of conical limits in \cite{ABLR},
and in full generality in \cite{Kelly-Schmitt}. It could equivalently be
stated as
\begin{quotation}
  $\Phi$-limits commute in \V with colimits that have $\Phi$-continuous 
  weights
\end{quotation}
or
\begin{quotation}
  All $\Phi$-continuous weights are $\Phi$-flat.
\end{quotation}

It is this condition 
which allows us to ``decompose'' colimits, in analogy with the finite/filtered
decomposition, where now $\Phi$-colimits play the role of ``finite'', and 
$\Phi$-flat colimits play the role of ``filtered''. This can be done
thanks to the following, which is a restatement of parts 
of \cite[Theorems~8.9,~8.11]{Kelly-Schmitt}:

\begin{theorem}\label{thm:KS}
The following condition on the class $\Phi$ of weights are equivalent:
\begin{enumerate}
\item If $F:\A\to\V$ is $\Phi$-continuous then so is $F*-:[\A\op,\V]\to\V$
(Axiom~A)
\item The category $\PhiCts(\A,\V)$ of $\Phi$-continuous presheaves is the free
completion of $\A\op$ under $\Phi$-flat colimits
\item Any presheaf $F:\A\to\V$ is a $\Phi$-flat colimit of presheaves in $\Phi$
\end{enumerate}
where in each case \A is allowed to be any small \V-category,
$\Phi$-complete in the first two cases.
\end{theorem}

\begin{proposition}\label{prop:gen-AxA}
Let \A and \B be small \V-categories with $\Phi$-limits, and 
$G:\A\to\B$ an arbitrary \V-functor. If $M:\A\to\V$ preserves 
$\Phi$-limits then so does $\Lan_G M:\B\to\V$.
\end{proposition}

\proof
There is a functor $\B(G,1):\B\to[\A\op,\V]$ sending $B\in\B$ to 
$\B(G-,B):\A\op\to\V$ which in turn sends an object $A\in\A$ to 
the hom-object $\B(GA,B)$. This functor preserves all existing limits,
and its composite with $\Lan_Y M:[\A\op,\V]\to\V$  is $\Lan_G M$.
\endproof

We now turn to the requirements on \K. It is possible to define
models in any \V-category with $\Phi$-limits, but in order to develop
the theory, somewhat more is required. We shall consider two levels 
of generality (more precise conditions will be given later).

\setcounter{axiom}{0}
\renewcommand{\theaxiom}{B\arabic{axiom}.}
\begin{axiom}
\K is locally $\Phi$-presentable: this is equivalent to saying that 
\K itself has the form $\PhiCts(\T,\V)$ for some small \V-category \T with 
$\Phi$-limits. It follows that $\K$ is reflective in $[\T,\V]$, and so is
complete and cocomplete.
\end{axiom}
\begin{axiom} \K has $\Phi$-limits, and the inclusion 
$y:\K\to\PK$ has a $\Phi$-continuous left adjoint.
\end{axiom}

By Theorem~\ref{thm:KS}, the Axiom~B1 is equivalent to saying that 
\K is the free completion under $\Phi$-flat colimits of a small 
$\Phi$-cocomplete \V-category. Axiom~B2 implies in particular that
\K is cocomplete; note that \PK has $\Phi$-limits by 
\cite[Proposition~4.3]{small} and Axiom~A. Axiom~B2 is a strong 
exactness condition, related to lex-totality \cite{yoneda}: when \V is \Set, and $\Phi$ consists of the finite
limits, it holds in any Grothendieck topos.

\begin{remark}
Both axioms imply that $\Phi$-limits commute in \K with $\Phi$-flat colimits,
since this is true in \V, and so in both $[\T,\V]$ and $\PK$, since the 
the limits and colimits are computed pointwise there. Now $\PhiCts(\T,\V)$ 
is closed in $[\T,\V]$ under limits and $\Phi$-flat colimits, so the
desired commutativity remains true there. In the case of Axiom~B2, both 
$\Phi$-limits and arbitrary colimits may be computed in \K by passing to 
\PK (where they commute) and then reflecting back into \K.
\end{remark}

\begin{proposition}
A presheaf category $\K=[\C,\V]$ satisfies Axioms~B1 and~B2.
\end{proposition}

\proof
Axiom~B1 is easy: if \T is the free \V-category with $\Phi$-limits on \C,
then $\PhiCts(\T,\V)\simeq[\C,\V]$.

As for Axiom~B2, since $[\C,\V]$ is cocomplete, the Yoneda functor 
$Y:[\C,\V]\to\P[\C,\V]$ certainly has a left adjoint $L\dashv Y$.
Explicitly, for a small presheaf $G:[\C,\V]\op\to\V$, the reflection
$LG\in[\C,\V]$ is the functor sending $C\in\C$ to $G(\C(C,-))$.

Writing $\ev_C$ for the functor $[\C,\V]$ given by evaluation at $C$,
and $yC$ for the representable $\C(C,-)$, we may therefore characterize
$L$ by the isomorphisms $\ev_C L\cong \ev_{yC}$, natural in $C$.

Let $F:\D\to\V$ be in $\Phi$, and $S:\D\to\P[\C\op,\V]$.
We must show that $L$ preserves the limit $\{F,S\}$; but this 
will be true if and only if $\ev_C L$ preserves the limit for each
$C\in\C$. This we show as follows:
\begin{align*}
\ev_C L\{F,S\} &\cong \ev_{yC}\{F,S\} \\
                  &\cong \{F,\ev_{yC}S\} \\
                  &\cong \{F,\ev_C LS\}
\end{align*}
using the fact that limits are preserved by evaluation functors.
\endproof

\section{Sifted colimits and locally \sfly presentable categories}\label{sect:sifted}

A key notion will be that of {\em sifted colimit} \cite{sifted, tamisant}.
A small category \D is 
said to be {\em sifted} if \D-colimits commute in \Set with finite products;
equivalently, if \D is non-empty and for all objects $A,B\in\D$, the category
of cospans from $A$ to $B$ is connected. Reflexive coequalizers and
filtered colimits are both sifted, but Ad\'amek \cite{what-is-sifted} has given an example of 
a category with reflexive coequalizers and filtered colimits but not all
sifted colimits.

Much of the theory of filtered colimits, involving finitely presentable
objects and locally finitely presentable categories, has analogues 
involving sifted colimits. This was developed in \cite{sifted}, and 
put into a more general setting in \cite{ABLR}.

For example, a functor $F:\A\to\Set$ is said to be 
{\em sifted-flat} \cite{sifted} if the left Kan extension $\Lan_Y F:[\A\op,\Set]\to\Set$ 
preserves finite products. Since $\Lan_Y F$ is also the functor 
$F*-:[\A\op,\Set]\to\Set$ calculating the $F$-weighted colimit of a 
functor $\A\op\to\Set$, to say that $F$ is sifted-flat is equivalently 
to say that $F$-weighted colimits commute in \Set with finite products.
In particular, if \A is sifted then $\Delta 1:\A\op\to\Set$ is 
sifted-flat. The sifted-flat functors were characterized in
\cite[Theorem~2.6]{sifted} as the presheaves which are sifted colimits
of representables.

An object $X$ of a category \K with sifted colimits is called 
{\em \sfly presentable} \cite{sifted} if the representable functor $\K(X,-):\K\to\Set$ preserves sifted colimits. 

In the following theorem,
the equivalence between (i) and (ii) is a special case of the general 
characterization
of free completions under colimits \cite[Proposition~5.62]{Kelly-book},
while the equivalence between (ii) and (iii) is \cite[Theorem~3.10]{sifted}, 
but is essentially already in \cite[Proposition~5.52]{Diers-density}.

\begin{theorem}\label{thm:variety}
For a category \K the following conditions are equivalent:
\begin{enumerate}[(i)]
\item \K is cocomplete and is the free completion under sifted colimits
of a small category \G
\item \K is cocomplete and has a small full subcategory \G consisting of
\sfly presentable objects, such that every object of \K is
a sifted colimit of objects in \G
\item \K is equivalent to the category $\FPP(\G\op,\Set)$ of 
finite-product-preserving functors from a small category \G with finite 
coproducts to \Set.
\end{enumerate}
The \sfly presentable objects will be the closure under retracts
of the category \G in each case.
\end{theorem}

Such a category \K is called \lsfp. The \lsfp categories are the 
(possibly multisorted) varieties.

We saw above that reflexive coequalizers and filtered colimits are not enough to guarantee
all sifted colimits. On the other hand, the following proposition shows that {\em preservation} of reflexive 
coequalizers and filtered colimits is enough to guarantee preservation of
sifted colimits, provided that all colimits actually exist. The history of
this result is slightly complicated. The 
first-named author knew the result and its proof from soon after the time of the first papers on sifted colimits, but did not know that it was regarded
as an important open problem, and did not publish it. An analogue in the context of quasicategories was
recently proved by Joyal \cite{Joyal-quasicats}, and inspired by this, the result itself was 
proved by Ad\'amek \cite{what-is-sifted}.

\begin{proposition}\label{prop:funny}
Let \K be a cocomplete category, and $F:\K\to\LL$ a functor. Then $F$
preserves sifted colimits if and only if it preserves reflexive coequalizers
and filtered colimits.
\end{proposition}

\proof
Since reflexive coequalizers and filtered colimits are both sifted colimits,
one direction is immediate. For the converse, suppose that $F$ preserves
reflexive coequalizers and filtered colimits.
Let \D be a sifted category, and $S:\D\to\K$ a diagram. We must show
that $F$ preserves the colimit of $S$. 

Let $\Fam\D$ be the free completion of \D under finite coproducts, and
$J:\D\op\to(\Fam\D)\op$ the canonical inclusion. Let $G=\Lan_J\Delta1$ be 
the left Kan extension of the terminal functor $\Delta1:\D\op\to\Set$ along
$J$. Since \D is sifted, $\Delta 1:\D\op\to\Set$ is sifted-flat,
hence so by Proposition~\ref{prop:flat-saturated} is its left Kan 
extension $G$, and so $G$ is a sifted colimit of representables.
But $\Fam\D$ has finite coproducts, and so sifted colimits can be constructed
using reflexive coequalizers and filtered colimits by 
\cite[Example~2.3(2)]{sifted}; thus $G$ can be 
constructed from the representables using reflexive coequalizers and filtered
colimits. We conclude that any functor preserving reflexive coequalizers
and filtered colimits also preserves $G$-weighted colimits.

For a functor $R$ with domain $\Fam\D$, we have (see Proposition~\ref{prop:Lan-weight}) canonical isomorphisms
$$G*R=\Lan_J\Delta1*R\cong \Delta1*RJ\cong\colim(RJ)$$
with all terms existing if any one of them does.
If $R:\Fam\D\to\K$ is the finite-coproduct-preserving functor extending
$S$, then we get $G*R\cong\colim(RJ)\cong\colim(S)$.
Since $F$ preserves reflexive coequalizers and filtered colimits, it preserves 
$G$-weighted colimits, and now 
$$F\colim(S)\cong F(G*R)\cong G*FR\cong \Delta1*FRJ\cong\colim(FRJ)\cong\colim(FS).$$
\endproof

\section{Possible choices for the class of limits}
\label{sect:examples}

\subsection{Finite limits}\label{sect:amiens}

Fix a symmetric monoidal closed category $\V=(\V_0,\ot,I)$ with underlying
category $\V_0$ complete and complete. As usual, an object $x\in\V_0$ is called 
finitely presentable if the representable functor $\V_0(x,-):\V_0\to\Set$ 
preserves filtered colimits.

Kelly \cite{Kelly-amiens} defines \V to be {\em locally finitely
presentable as a closed category} if $\V_0$ is locally finitely presentable
in the usual sense, and the finitely presentable objects are closed under
the monoidal structure: the unit $I$ is finitely presentable, and the tensor
product of any two finitely presentable objects is finitely presentable.

\begin{remark}
In fact all the key results of \cite{Kelly-amiens} remain true if we
drop the assumption that $I$ is finitely presentable: see 
Remark~\ref{rmk:amiens2} below. What is lost is the fact that finite
presentability in \V is the same as finite presentability in $\V_0$: one
only knows that every finitely presentable object in $\V_0$ is finitely
presentable in \V, not the converse. This possible generalization
seems to be of limited interest ---
we know of no important new examples --- but we mention it here, because
a similar generalization will be important 
when we move from the locally finitely presentable case to the \lsfp case. 
\end{remark}


This now gives a good notion of finite limit in the \V-enriched sense. 
First of all an object $X$ of a cocomplete \V-category \K is
said to be finitely presentable if the hom-functor $\K(X,-):\K\to\V$ preserves
filtered colimits. 

For an object $X$ of a cocomplete \V-category \K, there is in general no 
relation between the property of being finitely presentable in \K and the 
property of being finitely presentable in the underlying ordinary category
$\K_0$ of \K. But if \V is locally finitely presentable as a closed category,
then these two notions agree for $\K=\V$ (and more generally for any 
locally finitely presentable \V-category \K).

The {\em finite} limits are now those in the saturation of the class of
finite conical limits and \Vf-powers, where \Vf-powers are powers (cotensors)
by finitely presentable objects of \V.
We now take these finite limits to be our class $\Phi$. 
The fact that Axiom~A holds is Theorem~6.12 of 
\cite{Kelly-amiens}. Furthermore, Axiom~B will hold if \K is any locally 
finitely presentable \V-category, in the sense of \cite{Kelly-amiens}; in other 
words, if \K is a full reflective subcategory of a presheaf category
$[\C,\V]$ which is closed in $[\C,\V]$ under filtered colimits; or,
equivalently, if \K is the category $\Lex(\T,\V)$ of models of a small
\V-category \T with finite limits.

Many examples of locally finitely presentable \V were given in 
\cite{Kelly-amiens}; they include the closed categories of 
sets, pointed sets, abelian groups, modules over a commutative ring,
chain complexes, categories groupoids, and simplicial sets.



Most of the general results we shall prove about monads and theories for
$\K=\V$ 
were obtained in this case by Kelly in \cite{Kelly-amiens}, but for the 
part involving monads on, and single-sorted theories in, \V see
\cite{Power-EnrichedLawvere}. The treatment of monads and single-sorted
theories for more general \K (still with $\Phi$ the class of finite
limits) appeared in \cite{Power-Nishizawa}.

\subsection{Finite products}\label{sect:fp}

In \cite{borceux-day-univ-alg} the notion of $\pi$-category was introduced 
as a suitable setting for universal algebra. Explicitly, this is 
a complete and cocomplete symmetric monoidal category \V, such that 
Axiom~A holds for $\Phi$ the class of finite products, and furthermore,
the functors $-\t X:\V\to\V$ and $X\t-:\V\to\V$ preserve reflexive coequalizers
and filtered colimits for all objects $X$. By Proposition~\ref{prop:funny},
this condition on $-\t X$ and $X\t-$ is equivalent to saying that they
preserve {\em sifted colimits}.

\begin{proposition}
Let \K be a cocomplete category with finite products. The following 
conditions are equivalent:
\begin{enumerate}[(i)]
\item $X\t-:\K\to\K$ preserves sifted colimits for all $X$
\item $-\t X:\K\to\K$ preserves sifted colimits for all $X$
\item $\t:\K\t\K\to\K$ preserves sifted colimits
\item sifted colimits commute in \K with finite products
\end{enumerate}
\end{proposition}

In the setting of \cite{borceux-day-univ-alg}, the \V-category we are 
calling \K is always \V itself. Under these assumptions, the various
results we have considered were proved in \cite{borceux-day-univ-alg} (for the case $\K=\V$): left adjoints to algebraic functors,
reflectiveness of models, correspondence between monads and theories, and
so on.

Now finite products commute with sifted colimits in \Set, and more 
generally in any \lsfp category. 

If $\V_0$ is locally \sfly presentable then 
$-\t X$ and $X\t-$ will preserve sifted colimits; and now if Axiom~A holds,
\V will be a $\pi$-category.

In the case $\V=\Cat$, the $\Phi$-accessible monads on \Cat
will be the {\em strongly finitary 2-monads} of \cite{fpp}. These correspond
to (a finitary version of) the {\em discrete Lawvere theories} of 
\cite{Power:discrete-Lawvere}.

In the following section these sifted colimits will play a still 
more central role, and we shall see how they allow a more expressive notion of 
theory than that of \cite{borceux-day-univ-alg}, although with somewhat 
greater restrictions on \V.

\subsection{\Sfinite limits}\label{sect:lsfp}

In this section, which is one of the main original contributions of the paper,
we adapt the setting of \cite{Kelly-amiens} using sifted colimits in
place of filtered ones. The notion of \sfinite limit, introduced 
below, reduces to that of finite product in the case $\V=\Set$, but 
not in general. 




Suppose as usual that $\V=(\V_0,\ot,I)$ is a complete and cocomplete
symmetric monoidal closed category. This time we suppose that $\V_0$
is \lsfp, and so has the form $\FPP(\G\op,\Set)$ for a category \G with 
finite coproducts, which we may take to be the 
category of \sfly presentable objects of $\V_0$.

The directly analogous approach to \cite{Kelly-amiens} would be to suppose
that \G was closed under the monoidal structure; unfortunately this is not 
true for key examples such as $\V=\Gph$ (see Example~\ref{ex:graph}). 
We weaken the assumption slightly
by  supposing that \G is closed under the tensor product, but not that it 
contains the unit $I$. We then say that \V
is {\em locally \sfly presentable as a $\ot$-category}. 

Given such a \V, we can now develop the theory of locally strongly 
finitely presentable categories in the \V-enriched context. 
We say that an object $X$ of a cocomplete \V-category \K is {\em strongly
finitely presentable} if the hom-functor $\K(X,-):\K\to\V$ preserves sifted
colimits, and write \Ksf for the full subcategory of \K consisting of such
objects.

Just as in \cite{Kelly-amiens}, it is important to distinguish between
$X$ being \sfly presentable in \K, in the sense that 
$\K(X,-):\K\to\V$ preserves sifted colimits, and $X$ being strongly 
finitely presentable in $\K_0$, in the sense that $\K_0(X,-):\K_0\to\Set$
preserves such sifted colimits. (On the other hand, since \K and \V do 
have sifted colimits, to say that the \V-functor $\K(X,-):\K\to\V$ preserves 
sifted colimits is no different to saying that the underlying ordinary functor
$\K(X,-)_0:\K_0\to\V_0$ preserves sifted colimits.)
If $I\in\V_0$ were \sfp, then every 
\sfp object in \K would be \sfp in $\K_0$, but we are not assuming this.
What we do have is:

\begin{lemma}
If $X\in\V$ is \sfly presentable as an object of $\V_0$, in
the sense that $\V_0(X,-):\V_0\to\Set$ preserves sifted colimits, then it
is \sfly presentable in \V.  
\end{lemma}

\proof
Suppose that $X$ is \sfly presentable in $\V_0$. Then it is
a retract of a finite coproduct of objects in \G. If we now tensor $X$
by an arbitrary $G\in\G$, the resulting $G\ot X$ is again a retract of a 
finite coproduct of objects in \G, since \G is closed under tensoring. 
Thus $G\ot X$ is again \sfly presentable in $\V_0$, and so 
$\V_0(G\ot X,-)$ preserves sifted colimits. But this means that 
$\V_0(G,\V(X,-))$ preserves sifted colimits. Since the
$\V_0(G,-)$ preserve and jointly reflect sifted colimits, it follows that
$\V(X,-)$ preserves sifted colimits, and so that $X$ is \sfly
presentable in \V.
\endproof

The converse is false: we shall see in Example~\ref{ex:graph} 
below that if \V is the cartesian closed category \Gph of graphs, 
then the terminal object is \sfly presentable in \V but not
in $\V_0$.

As a further indication of the distinction between the properties 
of being \sfp in \V or $\V_0$, notice that although we had to {\em assume}
that the \sfly presentable objects of $\V_0$ were closed 
under tensoring, this is {\em automatic} for the \sfly presentable
objects of \V, since $\V(X,\V(Y,-))\cong\V(X\ot Y,-)$ and sifted-colimit-preserving functors are closed under composition.
 
An important technical result is:

\begin{proposition}
\Vsf is (equivalent to) a small \V-category.
\end{proposition}

\proof
Since $\V_0$ is \lsfp, it is 
also locally finitely presentable. Thus there is a regular cardinal $\alpha$
for which $I$ is $\alpha$-presentable; and $\V_0$ is still locally
$\alpha$-presentable. The $\alpha$-presentable objects 
of $\V_0$ are the $\alpha$-colimits of the finitely presentable ones, and
these are closed under tensoring, and by assumption they contain the unit 
object $I$. It now follows, just as in \cite[5.2, 5.3]{Kelly-amiens} that
an object is $\alpha$-presentable in \V if and only if it is 
$\alpha$-presentable in $\V_0$. An object of \Vsf is certainly 
$\alpha$-presentable in \V; thus $(\Vsf)_0$ is a full subcategory of 
$(\V_0)_\alpha$, which is small, and so \Vsf too is small.
\endproof

\begin{remark}\label{rmk:amiens2}
All the key results of \cite{Kelly-amiens} remain true if $\V_0$ is 
locally finitely presentable and the tensor product of two finitely 
presentable objects is finitely presentable. The reason for assuming that
the unit $I$ is finitely presentable is that then the notions of finite
presentability in \V and in $\V_0$ agree. But this is only used at one
point: in the proof of~7.1, in order to prove that the full subcategory
\Vf of finitely presentable objects in \V is small. However this can be
obtained alternatively as follows. As observed in \cite{vcat}, if 
$\V_0$ is locally finitely presentable, then \V is locally $\alpha$-presentable
as a closed category for some regular cardinal $\alpha$ --- by the same
argument that was used in the previous proposition. Then the full
subcategory $\V_\alpha$ of \V consisting of the $\alpha$-presentable objects
is small, by the same argument as in \cite{Kelly-amiens}; but \Vf is clearly
contained in $\V_\alpha$ and so is also small.
\end{remark}

Having fixed our monoidal category \V, we now turn to the class $\Phi$
of weights. We take for $\Phi$ the saturation of the class of finite products 
and \Vsf-powers (powers by objects of \Vsf).

\begin{proposition}\label{prop:pres}
If \K is cocomplete, the \sfly presentable objects of \K are
closed under $\Phi$-colimits.  
\end{proposition}

\proof
It suffices to show that they are closed under finite coproducts and under
\Vsf-copowers. If $X_1,\ldots,X_n$ are \sfly presentable, then
$$\K(X_1+\ldots+X_n,-)\cong\K(X_1,-)\t\ldots\t\K(X_n,-)$$
and each $\K(X_i,-)$ preserves sifted colimits since $X_i$ is strongly
finitely presentable, while 
finite products of sifted-colimit-preserving functors into \V still
preserve sifted colimits, since finite products commute with sifted 
colimits in \V. This proves that $X_1+\ldots+X_n$ is \sfly
presentable. 

Similarly if $X\in\K$ is \sfly presentable, and $G\in\V$ is
\sfly presentable, then
$$\K(G\cdot X,-) \cong \V(G,\K(X,-))$$
which preserves sifted colimits since $\K(X,-)$ and $\V(G,-)$ do, thus
$G\cdot X$ is \sfly presentable.
\endproof

\begin{corollary}
If $F:\A\op\to\V$ is in $\Phi$ then it is \sfly presentable
as an object of $[\A\op,\V]$.
\end{corollary}

\proof
This is immediate from the case $\K=[\A\op,\V]$ of the proposition,
given that representables are \sfly presentable, and $F$
is a $\Phi$-colimit of representables, via the Yoneda isomorphism
$F\cong F*Y$.
\endproof

The following theorem, adapted from \cite[Theorem~6.11]{Kelly-amiens}, 
implies in particular that Axiom~A holds:

\begin{theorem}\label{thm:sifted}
Let \T be a small \V-category with $\Phi$-limits. For a \V-functor
$F:\T\to\V$ the following are equivalent:
\begin{enumerate}[(1)]
\item $F$ is a sifted colimit of representables;
\item $F$ is $\Phi$-flat;
\item $F$ is $\Phi$-continuous.
\end{enumerate}
\end{theorem}

\proof
Representables are $\Phi$-flat, and sifted colimits commute in \V with
$\Phi$-limits, thus (1) implies (2). To see that (2) implies (3), observe
that if $F$ is $\Phi$-flat then $\Lan_Y F$ is $\Phi$-continuous, but
$Y$ is $\Phi$-continuous, hence so is $F=(\Lan_Y F)Y$.

So it remains to prove that (3) implies (1). Suppose then that $F$ is
$\Phi$-continuous. Consider the underlying ordinary functor 
$F_0:\T_0\to\V_0$, and the induced $\V_0(I,F_0):\T_0\to\Set$. Like any
\Set-valued functor, this is canonically a colimit of representables,
We form the category of elements \E and the induced $P:\E\to\T\op_0$. 
Explicitly, an object of \E consists of an object $T\in\T$ equipped with
a \V-natural $\T(T,-)\to F$. Then $\V_0(I,F_0)$ is the colimit of 
the composite 
$$\xymatrix{\E \ar[r]^P & \T\op_0 \ar[r]^Y & [\T_0,\Set] }$$
Since $F$ preserves finite products, so does $\V_0(I,F_0)$; it follows
that \E has finite coproducts, and so is sifted. Thus we have expressed
$\V_0(I,F_0)$ as a sifted colimit of representables.

The idea is to adapt this to obtain $F$ and not just $\V_0(I,F_0)$.
Consider now the composite
$$\xymatrix{ \E \ar[r]^P & \T\op_0 \ar[r]^{Y_0} & [\T,\V]_0 }$$
which sends an object $(T,x:\T(T,-)\to F)$ of \E to $\T(T,-)$.
We shall show that it has colimit $F$, and so that $F$ is a sifted colimit
of representables. There is an evident cocone $\gamma:Y_0 P\to\Delta F$, 
whose coprojection at
$(T,x:\T(T,-)\to F)$ is just $x$. We must show that this is a colimit.
It will be a colimit if and only if it gives a colimit after evaluating
at all $S\in\T$; in other words, if the \E-indexed diagram 
$$\xymatrix{ \T(T,S) \ar[r] & FS }$$
is a colimit (in \V, or equivalently in $\V_0$) for all $S$. Now the 
hom-functors $\V_0(G,-):\V_0\to\Set$ for $G\in\G$ jointly reflect colimits,
since \G is a strong generator, and they preserve sifted colimits, by
assumption on \G. Thus we are reduced to showing that the \E-indexed diagram
$$\xymatrix{ \V_0(G,\T(T,S)) \ar[r] & \V_0(G,FS) }$$
is a colimit in \Set, for all $S\in\T$ and all $G\in\G$. But by the universal
property of powers, and the fact that $F$ preserves \G-powers, this
diagram is equivalently
$$\xymatrix{ \V_0(I,\T(T,G\pitchfork S)) \ar[r] & \V_0(I,F(G\pitchfork S)) }$$
and this is a colimit since 
$$\xymatrix{ \V_0(I,\T(T,-)) \ar[r] & \V_0(I,F) }$$
is one. 
\endproof

We now turn to examples of \V which are \lsfp as $\ot$-categories. 
An important special case is where $\V$ is
a presheaf category, equipped with the cartesian closed structure. First
we prove:

\begin{lemma}
In a presheaf category $[\C\op,\Set]$, an object is a retract of a finite
coproduct of representables if and only if it is a finite coproduct of
retracts of representables. 
\end{lemma}

\proof
If $R_i$ is a retract of the representable $yC_i$ for $i=1,\ldots,n$ then
$\sum_i R_i$ is a retract of $\sum_i yC_i$. Thus a finite coproduct of 
retracts of representables is a retract of a finite coproduct of representables.

Conversely, suppose that $R$ is a retract of a finite coproduct 
$\sum_i yC_i$ of representables. By extensivity of $[\C\op,\Set]$,
we can write the inclusion $R\to\sum_i yC_i$ as a coproduct 
$\sum_i R_i\to \sum_i yC_i$ where each $R_i\to yC_i$ is the inclusion 
of a retract.
\endproof

\begin{corollary}\label{cor:retract1}
If idempotents split in \C, then the retracts of finite coproducts of
representables are just the finite coproducts of representables.
\end{corollary}

\begin{corollary}\label{cor:retract2}
A retract of a finite coproduct of retracts of finite coproducts of 
representables is a retract of a finite coproduct of representables. 
\end{corollary}

\begin{proposition}
If $\V=[\C\op,\Set]$, equipped with the cartesian closed structure, then
\V is locally \sfly presentable as a $\ot$-category if and 
only if the product of any two representables
is a  retract of a finite coproduct of representables.
\end{proposition}

\proof
In this case $(\V_0)_{sf}$ is obtained from \C by freely adjoining finite 
coproducts
and then splitting idempotents; in other words it may be identified with 
the full subcategory of $[\C\op,\Set]$ consisting of the retracts of finite
coproducts of representables. Then \V will be locally \sfly 
presentable as a $\ot$-category if and only if this subcategory is closed
under binary products. 

This certainly implies that the product of any
two representables is in the subcategory, but in fact it is equivalent:
if $R$ and $S$ are retracts of finite coproducts $\sum_i C_i$ and $\sum_j D_j$
of representables, then $R\times S$ is a retract of the finite sum 
$\sum_i\sum_j C_i\times D_j$ of products of representables. If each 
$C_i\times D_j$ is a retract of a finite coproduct of representables, then
so is $R\times S$ by Corollary~\ref{cor:retract2}.
\endproof

\begin{corollary}
$[\C\op,\Set]$ is \lsfp as a $\ot$-category if \C has binary products.  
\end{corollary}

\begin{example}\label{ex:graph}
The cartesian closed category \Gph of (directed) graphs is \lsfp as
a $\ot$-category. There are two representables: the free-living edge $E$,
and the free-living vertex $V$. It is easy to check that $V\t V\cong V$,
$V\t E\cong E\t V\cong V+V$, and $E\t E\cong V+E+V$, so that the product 
of any two representables is a finite coproduct of representables. On
the other hand, the terminal graph $1$ (which is of course the unit for 
the cartesian monoidal structure) is not a finite coproduct of representables,
and so \Gph is not \lsfp as a {\em closed} category. 

We can also see this more
directly by exhibiting a sifted colimit in \Gph which $\Gph(1,-):\Gph\to\Set$
fails to preserve. Now $\Gph(1,-)$ sends a graph $G$ to a vertex $v\in G$
with a chosen loop. Consider the parallel pair 
$$\xymatrix{
E+V \ar@<1ex>[r]^-{f} \ar@<-1ex>[r]_-{g} & E  }$$
where $f$ and $g$ act as the identity on the $E$ component, while on the 
$V$ component $f$ picks out the source of the edge in $E$, and $g$ picks
out the target. The coequalizer is formed by identifying the source and target, which now gives a loop: in other words, the coequalizer is $1$.
This pair is clearly reflexive, and so its coequalizer is a sifted colimit. 
But it is not preserved by $\Gph(1,-)$, since $\Gph(1,1)$ has one element, 
while $\Gph(1,E)$ is empty.
\end{example}

\begin{example}
The cartesian closed category \RGph of reflexive graphs is not \lsfp as a 
$\ot$-category. In this case there are two representables: the terminal graph 
$1$ with  a single loop, and the reflexive graph $E$ generated by a single 
edge. It follows
that any \sfly presentable reflexive graph can only have these
graphs as its connected components. But if we take the product of two copies
of $E$ we obtain the reflexive graph
$$
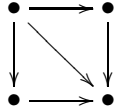
\xymatrix{
\bullet \ar[r] \ar[dr] \ar[d] & \bullet \ar[d] \\
\bullet \ar[r] & \bullet}$$
which, by Corollary~\ref{cor:retract1}, is not \sfly presentable.
(This time the terminal object is not just \sfly presentable, 
but representable.)
\end{example}

\begin{example}
Let \II be a skeletal category of finite sets and injections, and
\V the cartesian closed category $[\II,\Set]$. This is \lsfp as a 
$\ot$-category: given representables $\II(m,-)$ and $\II(n,-)$ we 
have the formula
$$\II(m,-)\t\II(n,-) \cong 
\sum^{m+n}_{k=\max(m,n)}\textstyle\binom{m}{m+n-k}\t\binom{n}{m+n-k}\cdot\II(k,-)$$
exhibiting $\II(m,-)\t\II(n,-)$ as a finite coproduct of representables.

The category [\II,Set] is one of models giving a denotational
semantics for the $\pi$-calculus \cite{Stark}. Since all operations
used there have \sfly presentable arities, our
formalism might be useful there.
\end{example}

It is now possible to develop a theory of locally \sfly 
presentable \V-categories as one might expect. We say that a \V-category
\K is \lsfp if
\begin{enumerate}[(i)]
\item \K is cocomplete
\item \K has a small full subcategory \G consisting of \sfly
 presentable objects
\item every object of \K is a sifted colimit of objects in \G
\end{enumerate}

We shall provide various characterizations below; in the meantime we
prove:

\begin{theorem}\label{thm:lsfp1}
Any \lsfp~\V-category \K is equivalent to the category of 
$\Phi$-continuous \V-functors from \T to \V for some small $\Phi$-complete
\V-category \T.
\end{theorem}

\proof
Let \GG be the closure of \G in \K under $\Phi$-colimits; that is,
under finite coproducts and \Vsf-copowers. Since \Vsf is small, 
it follows that \GG is still small.
By assumption it has $\Phi$-colimits, so $\T=\GG\op$ is a small
\V-category with $\Phi$-limits. The inclusion of $\GG\to\K$ induces
a \V-functor $W:\K\to[\GG\op,\V]=[\T,\V]$. Since $\Phi$-colimits of
\sfp objects are \sfp, the objects of \GG are \sfp, and so the 
functor $W:\K\to[\GG\op,\V]$ preserves sifted colimits.

Since \G is contained in \Ksf
it follows by  Proposition~\ref{prop:pres} that \GG is so too, and so
that $W$ actually lands in the category $\PhiCts(\T,\V)$ of $\Phi$-continuous
functors. Since \G is dense \cite[Theorem~5.35]{Kelly-book} it follows 
that \GG is dense and so that $W$ is fully faithful 
\cite[Theorem~5.13]{Kelly-book}. It remains to show 
that $\PhiCts(\T,\V)$ is the image
of $W$; in other words that every $\Phi$-continuous $F:\T\to\V$ has
the form $WA$ for some $A\in\K$. 

Suppose then that $F:\GG\op\to\V$ is $\Phi$-continuous;
by Theorem~\ref{thm:sifted} it is a sifted colimit of representables,
say $F=\colim_i yC_i$. Write $J:\GG\to\K$ for the inclusion. 
Since \K is cocomplete, we may form the colimit 
$F*J=(\colim_i yC_i)*J\cong\colim_i(yC_i*J)\cong\colim_i(JC_i)$
which will be preserved by $W$, since $W$ preserves sifted colimits.
Then $W(F*J)\cong F*WJ\cong F*Y\cong F$, and so $F$ is indeed in the image
of $W$.
\endproof

We shall see below that the converse is also true, but we shall prove
this using some of the more general theory developed in Section~\ref{sect:general} below.

\subsection{Sound doctrines}

We now show the case of finite limits and of \sfinite limits
fit into a common framework. The treatment follows that of 
Section~\ref{sect:lsfp} very closely, so we leave out some of the details.

Recall from \cite{ABLR} the notion of sound doctrine. In that paper a
doctrine consisted of a small collection \DD of small categories \D.
Then a category \C was said to be \DD-filtered if \C-colimits commute
in \Set with \DD-limits. It follows that for any $\D\in\DD$ and any
diagram $S:\D\op\to\C$, the category of cocones is connected. If conversely,
the connectedness of the category of cocones of a diagram $\D\op\to\C$ 
implies that \C is \DD-filtered, then the doctrine \DD is said to be 
{\em sound}.

The first main theorem about such doctrines \cite[Theorem~2.4]{ABLR} 
includes in particular:

\begin{theorem}
For a sound doctrine \DD and a functor $F:\A\to\Set$ with \A small,
the following are equivalent:
\begin{enumerate}
\item $\Lan_Y F$ is \DD-continuous
\item $F$ is a \DD-filtered colimit of representables
\end{enumerate}
and if \A has \DD-limits then these are further equivalent to 
\begin{enumerate}
\addtocounter{enumi}{2}
\item $F$ is \DD-continuous.
\end{enumerate}
\end{theorem}

We recover the setting of Sections~\ref{sect:amiens}
and~\ref{sect:lsfp} by taking for \DD, respectively, all finite
categories and all finite discrete categories.

A category is said to be {\em locally \DD-presentable} \cite{ABLR} if it is equivalent
to the category of \DD-continuous functors from \A to \Set, for a small
category \A with \DD-limits. There are analogues for all the main results
about locally finitely presentable categories, and these results are then 
recovered on taking \DD to be the finite categories.

Suppose that \DD is sound, and let $\V_0$ be locally \DD-presentable,
with the subcategory $(\V_0)_\DD$ of \DD-presentable objects closed 
under the tensor product in \V. We then say that \V is {\em locally
\DD-presentable as a $\ot$-category}. 

Let $\V_\DD$ be the full sub-\V-category 
of \V consisting of the \DD-presentable objects {\em of \V}: those objects
$G\in\V$ for which $[G,-]:\V\to\V$ preserves \DD-filtered colimits. 

\begin{lemma}
If $G$ is \DD-presentable in $\V_0$ then it is \DD-presentable in \V.
\end{lemma}

\proof
We know that $\V_0$ is (equivalent to) the category of \DD-continuous 
functors from \T to \Set for some small category \T with \DD-limits. 
If \DD-limits do not already include the splittings of idempotents then
we may split the idempotents of \T without changing the \DD-continuous
functors. A \DD-presentable object of $\V_0$ is then a representable
functor $yT=\T(T,-)$. We must show that the internal hom $[yT,-]:\V\to\V$
preserves \DD-filtered colimits, or equivalently that the underlying
functor $[yT,-]_0:\V_0\to\V_0$ does. But the $\V_0(yS,-):\V_0\to\Set$
preserve and detect \DD-filtered colimits, so it suffices to show that 
the $\V_0(yS,[yT,-])$ preserve \DD-filtered colimits. Finally 
$\V_0(yS,[yT,-])\cong\V_0(yS\ot yT,-)$, but $\V_0(yS\ot yT,-)$ preserves
\DD-filtered colimits since $yS\ot yT$ is \DD-presentable
by assumption.
\endproof

\begin{proposition}
$\V_\DD$ is (equivalent to) a small \V-category.
\end{proposition}

\proof
Any locally \DD-presentable category is locally $\alpha$-presentable for
some regular cardinal $\alpha$ (see \cite[Theorem~5.5]{ABLR}): we may take $\alpha$ to be larger than
the cardinality of any $\D\in\DD$ and such that the unit $I$ is 
$\alpha$-presentable. It then follows that \V is locally $\alpha$-presentable
as a closed category, and so that $\V_\alpha$ is small; but $\V_\DD$ is
contained in $\V_\alpha$.
\endproof

Let $\Phi$ be the saturation of the 
class of all (conical) \DD-limits and $\V_\DD$-powers. These limits certainly 
commute in \V with \DD-filtered colimits.

\begin{proposition}
If \K is a cocomplete \V-category, the $\Phi$-presentable objects of \K
are closed under $\Phi$-colimits. In particular, if $F:\A\op\to\V$
is in $\Phi$, then it is $\Phi$-presentable in $[\A\op,\V]$.
\end{proposition}

We now prove following analogue of Theorem~\ref{thm:sifted}, which
shows in particular that Axiom~A is satisfied. 

\begin{theorem}\label{thm:sound}
Let \T be a small \V-category with $\Phi$-limits. For a \V-functor
$F:\T\to\V$ the following are equivalent:
\begin{enumerate}[(1)]
\item $F$ is a \DD-filtered colimit of representables;
\item $F$ is $\Phi$-flat;
\item $F$ is $\Phi$-continuous.
\end{enumerate}
\end{theorem}

\proof
Since \DD-filtered colimits commute in \V with $\Phi$-limits,
the $\Phi$-flat weights are closed under \DD-filtered colimits.
Since representables are certainly $\Phi$-flat, we deduce that (1)
implies (2). Of course (2)
implies (3) since the Yoneda embedding preserves all existing limits.
So it remains to show that (3) implies (1).

Suppose then that \T is a small \V-category with $\Phi$-limits and
that $F:\T\to\V$ is $\Phi$-continuous. Just as in the proof of 
Theorem~\ref{thm:sifted}, we consider the ordinary functor
$$\xymatrix{
\T_0 \ar[r]^{F_0} & \V_0 \ar[r]^{\V_0(I,-)} & \Set }$$
its category of elements \E and the the induced $P:\E\to\T\op_0$, and 
observe that $\V_0(I,F_0)$ is canonically the colimit of
$$\xymatrix{
\E \ar[r]^-P & \T\op_0 \ar[r]^-{Y} & [\T,\Set]. }$$ 
Since $F$ preserves \DD-limits, so does $\V_0(I,F_0):\T_0\to\Set$; it follows that \E is \DD-filtered. If the colimit of
$$\xymatrix{
\E \ar[r]^-P & \T\op_0 \ar[r]^-{Y_0} & [\T,\V]_0 }$$
is $F$, then $F$ will be a \DD-filtered colimit of representables, as required.

The verification goes exactly as in the proof of Theorem~\ref{thm:sifted}.
\endproof

\subsection{Finite connected limits}

This is the case where \DD consists of the finite connected categories.
A category \K is locally \DD-presentable if and only if it is locally finitely
presentable and furthermore its category of finitely presentable objects is itself
the free completion under finite colimits of a full subcategory. The 
\DD-presentable objects are those which are both finitely presentable and
connected.

We suppose that $\V_0$ is such a category and that the finitely presentable
connected objects are closed under tensoring.  This time we take $\Phi$ to be 
the saturation of the class of all finite connected conical limits
and all \G-powers. 

\begin{example}
The cartesian closed categories \Gph, \RGph, \Cat, \Gpd, \SSet, and \CGTop
of graphs, reflexive graphs, categories, groupoids, simplicial sets, and compactly generated spaces are all examples; so is any presheaf topos
(such as \Gph, \RGph, and \SSet), and so, of course, is \Set.
\end{example}

A monad will be $\Phi$-accessible if and only if it is finitary and
preserves coproducts.

\section{Many-sorted theories}\label{sect:general}

Let \V and $\Phi$ be given, satisfying Axiom~A. To start with,
we allow \K to be an arbitrary \V-category with $\Phi$-limits,
but before long we shall suppose that it satisfies Axiom~B1 or~B2.
The most important case is $\K=\V$, which satisfies both Axiom~B1
and~B2.

\subsection{Theories and models}

A small \V-category \T with $\Phi$-limits is called a 
{\em $\Phi$-theory in \V}, or just a theory when $\Phi$ and \V
are understood. A {\em model} of \T in \K is a $\Phi$-continuous
\V-functor from \T to \K. 

The \V-category of models of \T in \K is the full subcategory 
$\PhiCts(\T,\K)$ of the functor
category $[\T,\K]$ consisting of the models. When $\K=\V$, we 
write simply $\PhiMod(\T)$ for $\PhiCts(\T,\V)$.

\subsection{Left adjoints to algebraic functors}

A {\em morphism of theories} is a $\Phi$-continuous \V-functor
$G:\S\to\T$. Composition with $G$ induces a \V-functor
$G^*:\PhiMod(\T)\to\PhiMod(\S)$; such a \V-functor is called
$\Phi$-algebraic, or just algebraic.

Such algebraic functors have left adjoints: given a model $M:\S\to\V$
we may form the left Kan extension $\Lan_G M:\T\to\V$ of $M$ along $G$
and by Proposition~\ref{prop:gen-AxA} this is $\Phi$-continuous, and so 
is a model of \T; it is easy to see that it has the required universal 
property.

(In fact the existence of a left adjoint holds much more generally; the 
point here is that it can be constructed via left Kan extension.)

We now turn to the case of a general \K. Once again $G$ induces 
a \V-functor $G^*:\PhiCts(\T,\K)\to\PhiCts(\S,\K)$; such a $G$ might
be called ``$\Phi$-algebraic relative to \K''.

\begin{proposition}
Let \A and \B be small \V-categories with $\Phi$-limits, and $G:\A\to\B$
an arbitrary \V-functor. If \K satisfies Axiom~B2 and $M:\S\to\K$ is
$\Phi$-continuous, then $\Lan_G M:\T\to\K$ is also $\Phi$-continuous.
\end{proposition}

\proof
Let $Y:\K\to\PK$ be the Yoneda embedding and $L\dashv Y$ its 
$\Phi$-continuous left adjoint. Of course $Y$ preserves $\Phi$-limits (and
any other existing limits).

Since $L$ is cocontinuous, it preserves left Kan extensions, and so 
$\Lan_G M\cong \Lan_G LYM\cong L\Lan_G YM$. Now $L$ is $\Phi$-continuous
by assumption, thus it will suffice to show that $\Lan_G YM$ is. 

In other words, we can work with \PK rather than \K. But in \PK both the
left Kan extensions and the $\Phi$-limits are computed pointwise, and so we
actually need only consider the case $\K=\V$; this is 
Proposition~\ref{prop:gen-AxA}.
\endproof

Thus when Axiom~B2 holds we can once again construct left adjoints to 
algebraic functors by Kan extension. Once again, the existence of the 
left adjoint holds much more generally, certainly whenever \K is 
locally presentable.

We shall see below that left adjoints to algebraic functors  include in 
particular free models for single-sorted theories.

\subsection{Reflectiveness of models}

Let \FT be the free completion of \T under $\Phi$-limits. Since \T has
$\Phi$-limits, the canonical inclusion $J:\T\to\FT$ has a right adjoint $R$,
and the algebraic functor $R^*:\PhiMod(\T)\to\PhiMod(\FT)\simeq[\T,\V]$ 
has a left adjoint by the previous result. But this is just the full inclusion
$\PhiMod(\T)\to[\T,\V]$. Thus the models form a full reflective subcategory
of the functor category, and in particular $\PhiMod(\T,\V)$ 
is complete and cocomplete.

\begin{theorem}
$\PhiMod(\T)$ is reflective in $[\T,\V]$ and so is complete and cocomplete.
It is closed in $[\T,\V]$ under all limits and under $\Phi$-flat colimits.
\end{theorem}

In fact the models will be reflective much more generally (see 
\cite[Chapter~6]{Kelly-book}), but our framework gives a simple construction,
which can be computed in practice (provided that colimits in the base
category \V can be computed). 

Note also that whenever the models are reflective
we also have adjoints to algebraic functors:  given
a morphism of theories $G:\S\to\T$ and a model $M:\S\to\K$, first take
the left Kan extension $\Lan_G M:\T\to\V$ and then reflect into models.

More generally, $\PhiMod(\T,\K)$ will be reflective in $[\T,\K]$ if \K
satisfies Axiom~B2.

\subsection{Characterization}

In this section we characterize \V-categories of models in \V.

Let \M be a \V-category with $\Phi$-flat colimits. We define an object
$M\in\M$ to be {\em $\Phi$-presentable} if the hom-functor
$\M(M,-):\M\to\V$ preserves $\Phi$-flat colimits.

We define a \V-category \M to be {\em locally $\Phi$-presentable} if it
is cocomplete and has a small full subcategory \G consisting of 
$\Phi$-presentable objects such that every object of \M is a 
$\Phi$-flat colimit of objects in \G. 

It follows that \M is the free completion of \G under $\Phi$-flat colimits.

Let \GG be the closure of \G in \M under $\Phi$-colimits. This is 
a small dense full subcategory consisting of $\Phi$-presentable objects, and
it is $\Phi$-cocomplete by construction. The inclusion $J:\GG\to\M$ induces
a \V-functor $W:\M\to[\GG\op,\V]$ which is fully faithful since \GG is
dense. It has a left adjoint sending $F:\GG\op\to\V$ to the colimit
$F*J\in\M$. The composite $WJ$ is the Yoneda embedding.

Explicitly, $W$ sends an object $M\in\M$ to $\M(J-,M)$ which is 
$\Phi$-continuous since the inclusion $J:\GG\to\M$ is $\Phi$-cocontinuous. 
Thus $W$ lands in $\PhiMod(\GG\op)$. Furthermore, $W$ preserves
$\Phi$-flat colimits, since \GG consists of $\Phi$-presentable objects.

\begin{proposition}
The \V-functor $W:\M\to\PhiMod(\GG\op)$ is an equivalence.  
\end{proposition}

\proof
We already know that $W$ is fully faithful, so it will suffice to show
that it is essentially surjective on objects. Suppose then that 
$F:\GG\op\to\V$ is $\Phi$-continuous. By Axiom~A it is also $\Phi$-flat,
and so $W$ preserves $F$-weighted colimits. Now every presheaf is 
a colimit of representables, weighted by itself, and so we have
$$F\cong F*Y\cong F*WJ\cong W(F*J)$$
which completes the proof.
\endproof

Thus any locally $\Phi$-presentable \V-category is the category of 
models in \V for a small theory. Conversely, let \T be a theory
and consider $\PhiMod(\T)$. This is reflective in $[\T,\V]$ and so
cocomplete. The representables provide a small dense subcategory

Since $\Phi$-limits commute in \V with $\Phi$-flat colimits, the inclusion
$\PhiMod(\T)\to[\T,\V]$ preserves $\Phi$-flat colimits, which is equivalent
to saying that the representables are $\Phi$-presentable in 
$\PhiMod(\T)$. Finally if $F:\T\to\V$ is $\Phi$-continuous then it 
is $\Phi$-flat by Axiom~A, and so once again 
$F\cong F*Y\cong F*WJ\cong W(F*J)$ shows that $F$ is a $\Phi$-flat colimit
in $\PhiMod(\T)$ of representables. This proves that $\PhiMod(\T)$ is
locally $\Phi$-presentable, and so gives:

\begin{theorem}
A category \M is locally $\Phi$-presentable if and only if it is equivalent
to a category of $\Phi$-continuous \V-functors from \T to \V for a small
$\Phi$-complete \V-category \T.
\end{theorem}

As observed above, if $\Phi$ consists of the finite limits, then this was proved
in \cite{Kelly-amiens}; and, in the case $\V=\Set$, is of course due to 
Gabriel-Ulmer \cite{Gabriel-Ulmer}.
We now return to the special case of Section~\ref{sect:lsfp}.

\begin{theorem}\label{thm:lsfp}
Let \V be \lsfp as a $\ot$-category and let $\Phi$ be the saturation 
of the finite products and \Vsf-powers. For a \V-category \M the 
following are equivalent:
\begin{enumerate}[(i)]
\item \M is \lsfp
\item \M is locally $\Phi$-presentable
\item There is a small \V-category \T with finite products and 
\Vsf-powers for which \M is equivalent to the full subcategory 
of $[\T,\V]$ consisting of the \V-functors which preserve these limits.
\end{enumerate}
\end{theorem}

\proof
The equivalence of (ii) and (iii) is a special case of the previous theorem.
The fact that (i) implies (iii) is Theorem~\ref{thm:lsfp1}. To see that (iii)
implies (i), observe that a \M satisfying (iii) is cocomplete, since (iii) 
implies (ii), and now consider the small full subcategory consisting of the 
representables $\T\to\V$.
Certainly these are \sfly presentable. It remains to show
that every $\Phi$-continuous $F:\T\to\V$ is a sifted colimit of 
representables. But this is Theorem~\ref{thm:sifted}.
\endproof

This reduces to Theorem~\ref{thm:variety} in the case $\V=\Set$; see the discussion
before that theorem for the history of the result in that case.

\section{Lawvere theories}\label{sect:Lawvere}

In this section we turn to theories which can be thought of as
single-sorted, and see how to extend the classical correspondence between
such theories and monads. 

A major difference is that the resulting theories need not have all 
the limits under consideration. This was observed in
\cite{Power-Nishizawa} in the case where $\Phi$ consists of the finite
limits, and our approach is modelled on that of \cite{Power-Nishizawa}.
As a consequence, in this more general setting, a Lawvere theory need not
be a theory; we shall see, however, how a Lawvere theory does generate
a theory with the same models.

Let \A and \B be \V-categories with $\Phi$-flat colimits. We shall say that
a \V-functor $F:\A\to\B$ is {\em $\Phi$-accessible} if it preserves $\Phi$-flat
colimits. A monad on \A will be called $\Phi$-accessible
if its underlying endofunctor is so. We write $\Mnd_\Phi(\A)$ for the 
category of $\Phi$-accessible monads on \A.

\subsection{Lawvere $\Phi$-theories}

Let \K be a \V-category satisfying Axiom~B1.
Let $J:\K_\Phi\to\K$ be the full subcategory of \K consisting of the 
$\Phi$-presentable objects; then $\K_\Phi$ has $\Phi$-colimits and $J$
preserves them. Furthermore, \K is equivalent to the category of
$\Phi$-continuous \V-functors from $\K\op_\Phi$ to \V.

Let $T=(T,m,i)$ be a $\Phi$-accessible \V-monad on \K, write $\K^T$ for the
Eilenberg-Moore \V-category, with forgetful functor $U^T:\K^T\to\K$ and 
left adjoint $F^T\dashv U^T$.

We may factorize the composite $F^T J:\K_\Phi\to\K^T$ as an identity-on-objects
\V-functor $E:\K_\Phi\to\G$
followed by a fully faithful $H:\G\to\K^T$. The opposite of the resulting 
\V-category \G will become the Lawvere theory \LL corresponding to $T$.

Now $\K_\Phi$ has $\Phi$-colimits, preserved by $J$, while $F^T$ preserves
all colimits, so the composite $HE=F^T J:\V_\Phi\to\V^T$ preserves
$\Phi$-colimits. It follows that $E$ preserves $\Phi$-colimits, but it
does not follow that \G has all $\Phi$-colimits. It does have 
$\Phi$-colimits of diagrams in the image of $E$, but need not have 
$\Phi$-colimits in general. 

Note, however, that \G will have $F$-weighted colimits for any $F\in\Phi$ which is a weight for coproducts or copowers, since these involve only the objects of \G, and so the resulting diagrams will be in the image of the 
identity-on-object \V-functor $E$.

\begin{remark}
Thus in the special case of Section~\ref{sect:lsfp} where all weights
in $\Phi$ are of this type, \G will have $\Phi$-colimits. The same is true if
$\Phi$ consists of just the finite products.
\end{remark}

But in general, we make the following definition, given in 
\cite{Power-Nishizawa} for the case 
where $\Phi$ is the finite limits.

\begin{definition}
A Lawvere $\Phi$-theory in \V is an identity-on-object \V-functor
$E:\V\op_\Phi\to\LL$ which preserves $\Phi$-limits. A morphism of Lawvere
$\Phi$-theories is a commutative triangle of (identity-on-object) 
\V-functors; we write $\Law_\Phi(\K)$ for the resulting category of 
Lawvere $\Phi$-theories on \K.
\end{definition}

We cannot in general simply define a model to be a $\Phi$-continuous
\V-functor with domain \LL, since \LL may not have all $\Phi$-limits.
Instead we define the \V-category of models of \LL by the following
pullback in \VCat
$$\xymatrix{
\PhiCts(\LL,\V) \ar[r] \ar[d]_U & [\LL,\V] \ar[d]^{[E,\V]} \\
\K \ar[r]_{\K(J,1)} & [\K\op_\Phi,\V] }$$

\begin{remark}
As observed in \cite{JFP} in the case of finite limits, 
since \K is equivalent to the \V-category
of $\Phi$-continuous \V-functors from $\K\op_\Phi$ to \V, up to an 
equivalence, a model of \LL is just a \V-functor $M:\LL\to\V$ whose
restriction along $E$ is $\Phi$-continuous.
\end{remark}

\begin{remark}
If $\Phi$ consists only of finite products and/or powers, then $\Phi$-limits
in \LL are determined by those in $\K\op_\phi$, and so the restriction of
$M$ along $E$ is $\Phi$-continuous if and only if $M$ is $\Phi$-continuous.
\end{remark}

\begin{proposition}
$U$ is monadic via a $\Phi$-accessible monad.
\end{proposition}

\proof
$[\LL,\V]$ has all coequalizers,
$[E,\V]$ has both adjoints (given by left and right Kan extension), and 
$[E,\V]$ is also conservative because $E$ is the bijective on objects.
An easy application of Beck's theorem shows that $[E,\V]$ is monadic.
The pullback $U$ of $[E,\V]$ will still satisfy the conditions of Beck's
theorem, and so be monadic, provided that it has a left adjoint. 

Now $\K(J,1):\K\to[\K\op_\Phi,\V]$ is fully faithful and has a left adjoint,
$[E,\V]$ has a left adjoint, and the inclusion $\PhiCts(\LL,\V)\to[\LL,\V]$ 
has a left adjoint, thus $U$ does indeed have a left adjoint and so is
monadic. 

It remains to show that the monad is $\Phi$-accessible, or equivalently 
that $U$ preserves $\Phi$-flat colimits. But this follows because the
other three functors in the definition of $\PhiCts(\LL,\V)$ preserve
$\Phi$-flat colimits, and $\K(J,1)$ is fully faithful so also reflects 
them.
\endproof

Thus to every Lawvere theory we have associated a $\Phi$-accessible monad.
This gives the object-part of a functor $\mnd:\Law_\Phi(\K)\to\Mnd_\Phi(\K)$.

Conversely, for a $\Phi$-accessible monad $T$ on \K, the inclusion 
$H:\LL\op=\G\to\K^T$ induces a \V-functor $\K^T(H,1):\K^T\to[\LL,\V]$.

\begin{theorem}
The \V-functor $\K^T(H,1):\K^T\to[\LL,\V]$ restricts to an isomorphism
of \V-categories 
$$\K^T\simeq\PhiCts(\LL,\V).$$  
\end{theorem}

\proof
Composition with $E:\K\op_\Phi\to\G\op$ induces a \V-functor 
$$\xymatrix{[\G\op,\V] \ar[r]^{[E,\V]} & [\K\op_\Phi,\V] }$$
which has both adjoints, given by left and right Kan extension. 
Since $E$ is bijective on objects, $[E,\V]$ is conservative, and
now by the Beck theorem, it is monadic. Write $S$ for the induced monad.

Consider what happens when
we restrict the induced monad $S$ along the fully faithful 
$\K(J,1):\K\to[\K\op_\Phi,\V]$. We have 
\begin{align*}
(\Lan_E\K(J,X))E &\cong \int^{c\in\K_\Phi} \G(E-,Ec)\cdot \K(Jc,X) \\
                 &\cong \int^c \K^T(HE-,HEc)\cdot\K(Jc,X) 
\tag{because $H$ is fully faithful} \\
                 &= \int^c \K^T(F^T J-,F^T Jc)\cdot\K(Jc,X) \\
                 &\cong \int^c \K(J-,TJc)\cdot\K(Jc,X) 
\tag{by adjointness} \\
                 &\cong \K(J,X)*\bigl(\K(J,1)T\bigr) 
\tag{by the coend formula for weighted colimits} \\
                 &\cong \K(J,1)\bigl(\K(J,X)*T\bigr) 
\tag{because $\K(J,1)$ preserves $\Phi$-flat colimits and $\K(J,X)$
  is $\Phi$-flat} \\
                 &\cong \K(J,1)\bigl(TX\bigr) 
\tag{because $T$ is $\Phi$-accessible} \\
                 &= \K(J,TX)
\end{align*}
Thus the functor part of $S$
restricts to $T$. In fact the monad itself restricts, and so we 
conclude that a $T$-algebra is an $S$-algebra (an object of $[\G\op,\V]$)
whose underlying object (restriction along $E$) is in the image of 
$\K(J,1)$. But this is exactly the definition of models of \LL.

Similarly, a morphism of $T$-algebras is the same as a morphism of the
corresponding $S$-algebras, whence the result.
\endproof

We have associated a Lawvere $\Phi$-theory \LL to every 
$\Phi$-accessible \V-monad on \K. This process is clearly functorial, 
giving a functor $\th:\Mnd_\Phi(\K)\to\Law_\Phi(\K)$. 

\begin{theorem}
The functors $\mnd$ and $\th$ form an equivalence of categories
$\Mnd_\Phi(\K)\simeq\Law_\Phi(\K)$.
\end{theorem}

\proof
The previous theorem gives an 
isomorphism $\mnd\circ\th\cong1$. For the other isomorphism 
$\th\circ\mnd\cong1$, let $E:\K\op_\Phi\to\LL$ be a Lawvere theory,
and $T$ the induced monad $\mnd(\LL)$.
Then $\th(\mnd(\LL))$ can be obtained by factorizing 
$F^T J:\K_\Phi\to\PhiCts(\LL,\V)$ as an identity-on-object functor followed
by a fully faithful one. Now to form $F^TJc$, for $c\in\K_\Phi$, we
send $Jc$ to $\K(J,Jc):\K\op_\Phi\to\V$,
then to its left Kan extension $\Lan_E\K(J,Jc):\LL\to\V$, and then reflect this
into $\PhiCts(\LL,\V)$. But $\K(J,Jc)\cong\K_\Phi(-,c)$, whose left Kan 
extension along $E$ is $\LL(-,Ec)$, and this is already in $\PhiCts(\LL,\V)$.
But since $F^T J$ sends $c$ to $\LL(-,Ec)$, its identity-on-object/fully-faithful
factorization gives just \LL.
\endproof

We now turn to the $\Phi$-theory generated by a Lawvere $\Phi$-theory \LL.
Every representable functor $\LL(L,-):\LL\to\V$ is a model of \LL, and so 
we get a fully faithful embedding $Y:\LL\op\to\PhiCts(\LL,\V)$. Form
the closure of the representables in $\PhiCts(\LL,\V)$ under $\Phi$-colimits.
This gives fully faithful $K:\LL\to\T$ and $P:\T\to\PhiCts(\LL,\V)\op$.
Clearly $H$ preserves $\Phi$-limits, while $P$ preserves those $\Phi$-limits
in the image of $E$. 

Now \T is a small \V-category with $\Phi$-limits; that is, a $\Phi$-theory.
Furthermore, it has the same models as \LL. This is really a special case
of \cite[Proposition~6.23]{Kelly-book}, but we outline here the argument.

First of all the composite
$$\xymatrix{
\LL\op \ar[r]^{K} & \T\op \ar[r]^-{P} & \PhiCts(\LL,\V) }$$
is dense, and both $K$ and $P$ are fully faithful. It follows by 
\cite[Theorem~5.13]{Kelly-book} that both $K$ and $P$ are dense, and that 
$P\cong\Lan_K(PK)$. Since $P$ is dense, the induced functor
$\PhiCts(\LL,\V)(P,1):\PhiCts(\LL,\V)\to[\T,\V]$ is fully faithful. We must show
that its image is exactly the $\Phi$-continuous functors.

Now $P:\T\op\to\PhiCts(\LL,\V)$ is $\Phi$-cocontinuous, by construction
of $\T\op$, and so the induced $\PhiCts(\LL,\V)(P-,M):\T\to\V$ will be $\Phi$-continuous for all models
$M$. This proves that $\PhiCts(\LL,\V)(P,1)$ takes values among the 
$\Phi$-continuous functors. Conversely, let $G:\T\to\V$ be 
$\Phi$-continuous. Then $GKE$ is $\Phi$-continuous, and so $GK$ is
(isomorphic to) a model; and indeed the model can be calculated 
as $GK*P$. Thus
$$GK\cong \PhiCts(\LL,\V)(PK-,GK*P).$$
But $G$ and $\PhiCts(\LL,\V)(P-,GK*P)$ are both $\Phi$-continuous \V-functors
$\T\to\V$, and \T is the closure of \LL under $\Phi$-limits, thus 
$G$ will be isomorphic to $\PhiCts(\LL,\V)(P-,GK*P)$ provided their restrictions
along $K$ are isomorphic; but this is the previous displayed equation.




\bibliographystyle{plain}


\end{document}